\documentstyle{amsppt}
\voffset-10mm
\magnification1200
\pagewidth{130mm}
\pageheight{204mm}
\hfuzz=2.5pt\rightskip=0pt plus1pt
\binoppenalty=10000\relpenalty=10000\relax
\TagsOnRight
\loadbold
\nologo

\let\ge\geqslant
\let\epsilon\varepsilon
\redefine\d{\roman d}
\define\Log{\operatorname{Log}}
\topmatter
\title
Binomial sums related \\
to rational approximations to $\zeta(4)$
\endtitle
\author
Wadim Zudilin\footnotemark"$^\ddag$"\ {\rm(Moscow)}
\endauthor
\date
\hbox to100mm{\vbox{\hsize=100mm%
\centerline{E-print \tt math.CA/0311196}
\smallskip
\centerline{November 2003}
}}
\enddate
\address
\hbox to70mm{\vbox{\hsize=70mm%
\leftline{Moscow Lomonosov State University}
\leftline{Department of Mechanics and Mathematics}
\leftline{Vorobiovy Gory, GSP-2, Moscow 119992 RUSSIA}
\leftline{{\it URL\/}: \tt http://wain.mi.ras.ru/index.html}
}}
\endaddress
\email
{\tt wadim\@ips.ras.ru},
{\tt wzudilin\@mi.uni-koeln.de}
\endemail
\abstract
For the solution $\{u_n\}_{n=0}^\infty$
to the polynomial recursion
$(n+1)^5u_{n+1}-\allowmathbreak
3(2n+1)(3n^2+3n+1)(15n^2+15n+4)u_n
-3n^3(3n-1)(3n+1)u_{n-1}=0$, where $n=1,2,\dots$,
with the initial data $u_0=1$, $u_1=12$, we prove
that all $u_n$ are integers. The numbers $u_n$, $n=0,1,2,\dots$,
are denominators of rational approximations to~$\zeta(4)$
(see {\tt math.NT/0201024}). We use Andrews's generalization
of Whipple's transformation of a terminating ${}_7F_6(1)$-series
and the method from {\tt math.NT/0311114}.
\endabstract
\endtopmatter
\leftheadtext{W.~Zudilin}
\rightheadtext{Binomial sums related to $\zeta(4)$}
\footnotetext"$^\ddag$"{The work is supported
by an Alexander von Humboldt research fellowship.}
\document

Consider the following 3-term polynomial recursion:
$$
\align
&
(n+1)^5u_{n+1}
-3(2n+1)(3n^2+3n+1)(15n^2+15n+4)u_n
\\ &\qquad
-3n^3(3n-1)(3n+1)u_{n-1}=0
\qquad\text{for}\quad n\ge1,
\endalign
$$
and take the two linearly independent solutions
$\{u_n\}_{n=0}^\infty$ and $\{v_n\}_{n=0}^\infty$
determined by the inicial conditions
$u_0=1$, $u_1=12$ and $v_0=0$, $v_1=13$. In~\cite{Z1},
we give a hypergeometric interpretation of the sequence
$u_n\zeta(4)-v_n$, $n=0,1,2,\dots$, from which one obtains
the limit
$$
\lim_{n\to\infty}\frac{v_n}{u_n}=\zeta(4)=\frac{\pi^4}{90}
$$
and the representation
$$
\align
u_n
&=(-1)^{n+1}\sum_{l=0}^n\frac{\d}{\d l}
\Bigl(\frac n2-l\Bigr)\binom nl^4\binom{n+l}n^2\binom{2n-l}n^2
\\
&=(-1)^n\sum_{l=0}^n\Bigl(\frac n2-l\Bigr)
\binom nl^4\binom{n+l}n^2\binom{2n-l}n^2
\\ &\qquad\times
\biggl(\frac1{n/2-l}-6H_{n-l}+6H_l
-2H_{n+l}+2H_{2n-l}\biggr),
\tag1
\endalign
$$
where $H_l=\sum_{j=1}^lj^{-1}$ are harmonic numbers.
The integrality of all $u_n$ (conjectured in~\cite{Z1})
is not an immediate consequence of formula~\thetag{1}.
In the recent work~\cite{KR} C.~Krattenthaler and T.~Rivoal
prove (among several other useful theorems and beautiful
binomial identities) that
$$
u_n=\sum_{i,j}\binom ni^2\binom nj^2\binom{n+j}n
\binom{n+j-i}n\binom{2n-i}n,
\qquad n=0,1,2,\dots,
$$
from which one has the desired inclusions $u_n\in\Bbb Z$.
The main objective of the present note is to give a simpler
proof of the formula for the numbers $u_n$ as well as to indicate
some other representations that also show that all $u_n$ are
integers.

We use the standard notation
$$
{}_{r+1}F_r\biggl(\matrix\format&\,\c\\
a_0, & a_1, & \dots, & a_r \\ & b_1, & \dots, & b_r
\endmatrix\biggm|z\biggr)
=\sum_{l=0}^\infty\frac{(a_0)_l(a_1)_l\dotsb(a_r)_l}
{l!\,(b_1)_l\dotsb(b_r)_l}z^l
$$
for the generalized hypergeometric series;
the notation $(a)_l=a(a+1)\dotsb(a+l-1)$ for $l=1,2,\dots$
and $(a)_0=1$ stands for the Pochhammer symbol.

The following formula is due G.\,E.~Andrews.
Making the passage $q\to1$ in \cite{A, Theorem~4}
(see also \cite{Z2} for a related application of the identity)
we have: {\it for $s\ge1$ and $m$ a non-negative integer},
$$
\align
&
{}_{2s+3}F_{2s+2}\biggl(\matrix\format&\,\c\\
a, & 1+\frac12a, & b_1, & c_1, & b_2, & c_2, & \dots \\
& \frac12a, & 1+a-b_1, & 1+a-c_1, & 1+a-b_2, & 1+a-c_2, & \dots
\endmatrix
\\ &\qquad\qquad\qquad\qquad\qquad\qquad\qquad
\matrix\format&\,\c\\
\dots, & b_s, & c_s, & -m \\
\dots, & 1+a-b_s, & 1+a-c_s, & 1+a+m
\endmatrix\biggm|1\biggr)
\\ &\quad
=\frac{(1+a)_m(1+a-b_s-c_s)_m}{(1+a-b_s)_m(1+a-c_s)_m}
\sum_{l_1\ge0}\frac{(1+a-b_1-c_1)_{l_1}(b_2)_{l_1}(c_2)_{l_1}}
{l_1!\,(1+a-b_1)_{l_1}(1+a-c_1)_{l_1}}
\\ &\quad\qquad\times
\sum_{l_2\ge0}\frac{(1+a-b_2-c_2)_{l_2}(b_3)_{l_1+l_2}(c_3)_{l_1+l_2}}
{l_2!\,(1+a-b_2)_{l_1+l_2}(1+a-c_2)_{l_1+l_2}}
\dotsb
\\ &\quad\qquad\times
\sum_{l_{s-1}\ge0}
\frac{(1+a-b_{s-1}-c_{s-1})_{l_{s-1}}
(b_s)_{l_1+\dots+l_{s-1}}(c_s)_{l_1+\dots+l_{s-1}}}
{l_{s-1}!\,(1+a-b_{s-1})_{l_1+\dots+l_{s-1}}
(1+a-c_{s-1})_{l_1+\dots+l_{s-1}}}
\\ &\quad\qquad\qquad\times
\frac{(-m)_{l_1+\dots+l_{s-1}}}
{(b_s+c_s-a-m)_{l_1+\dots+l_{s-1}}}.
\endalign
$$

Taking $s=3$, $a=-n-2\epsilon$, $b_1=b_2=b_3=c_2=-n-\epsilon$,
$c_1=c_3=n-\epsilon+1$ and $m=n$, $i=l_1$, $j=l_1+l_2$,
we derive from Andrews's formula
$$
\align
&
\sum_{l=0}^n\frac{-\frac n2-\epsilon+l}{-\frac n2-\epsilon}
\cdot\frac{(-n-2\epsilon)_l}{(1)_l}
\cdot\frac{(-n)_l}{(1-2\epsilon)_l}
\cdot\biggl(\frac{(1+n-\epsilon)_l}{(-2n-\epsilon)_l}\biggr)^2
\cdot\biggl(\frac{(-n-\epsilon)_l}{(1-\epsilon)_l}\biggr)^4
\\ &\quad
=\frac{(1-n-2\epsilon)_n(-n)_n}{(1-\epsilon)_n(-2n-\epsilon)_n}
\sum_i\frac{(-n)_i(-n-\epsilon)_i^2}{i!\,(1-\epsilon)_i(-2n-\epsilon)_i}
\\ &\quad\qquad\times
\sum_j\frac{(1+n)_{j-i}(-n-\epsilon)_j(1+n-\epsilon)_j(-n)_j}
{(j-i)!\,(1-\epsilon)_j^2j!}.
\tag2
\endalign
$$
Using the trivial equality
$$
(1-n-2\epsilon)_n
=\frac{-\epsilon}{-\frac n2-\epsilon}\cdot(-n-2\epsilon)_n,
$$
we may rewrite \thetag{2} in the form
$$
\align
&
\frac1\epsilon\sum_{l=0}^nA_l(\epsilon)
\\ &\quad
=\frac1\epsilon\sum_{l=0}^n\Bigl(\frac n2+\epsilon-l\Bigr)
\cdot\frac{(-n-2\epsilon)_l}{(1)_l}
\cdot\frac{(-n)_l}{(1-2\epsilon)_l}
\cdot\biggl(\frac{(1+n-\epsilon)_l}{(-2n-\epsilon)_l}\biggr)^2
\cdot\biggl(\frac{(-n-\epsilon)_l}{(1-\epsilon)_l}\biggr)^4
\allowdisplaybreak &\quad
=\frac{(1-n-2\epsilon)_n(-n)_n}{(1-\epsilon)_n(-2n-\epsilon)_n}
\sum_i\frac{(-n)_i(-n-\epsilon)_i^2}{i!\,(1-\epsilon)_i(-2n-\epsilon)_i}
\\ &\quad\qquad\times
\sum_j\frac{(1+n)_{j-i}(-n-\epsilon)_j(1+n-\epsilon)_j(-n)_j}
{(j-i)!\,(1-\epsilon)_j^2j!}.
\tag3
\endalign
$$
Now, we tend $\epsilon$ to~$0$. On the right hand side of~\thetag{3}
we only need to plug $\epsilon=0$. To proceed with the left hand side,
we first note that $A_l(0)=-A_{n-l}(0)$ for all $l=0,1,\dots,n$,
hence
$$
\lim_{\epsilon\to0}\sum_{l=0}^nA_l(\epsilon)
=\sum_{l=0}^nA_l(0)=0
\tag4
$$
and we may apply the l'H\^opital rule:
$$
\align
&
\lim_{\epsilon\to0}\frac1\epsilon\sum_{l=0}^nA_l(\epsilon)
=\sum_{l=0}^n\frac{\partial A_l(\epsilon)}{\partial\epsilon}
\bigg|_{\epsilon=0}
=\sum_{l=0}^nA_l(0)\cdot\frac{\partial\Log A_l(\epsilon)}{\partial\epsilon}
\bigg|_{\epsilon=0}
\\ &\quad
=\sum_{l=0}^nA_l(0)
\cdot\biggl(\frac1{\frac n2-l}-2\sum_{j=1}^l\frac1{-n+j-1}
+2\sum_{j=1}^l\frac1j-2\sum_{j=1}^l\frac1{n+j}
\\ &\quad\qquad
+2\sum_{j=1}^l\frac1{-2n+j-1}
-4\sum_{j=1}^l\frac1{-n+j-1}+4\sum_{j=1}^l\frac1j\biggr)
\allowdisplaybreak &\quad
=\sum_{l=0}^nA_l(0)
\cdot\biggl(\frac1{\frac n2-l}+6(H_n-H_{n-l})+6H_l
-2(H_{n+l}-H_n)-2(H_{2n}-H_{2n-l})\biggr)
\allowdisplaybreak &\quad
=\sum_{l=0}^nA_l(0)
\cdot\biggl(\frac1{\frac n2-l}-6H_{n-l}+6H_l
-2H_{n+l}+2H_{2n-l}\biggr),
\endalign
$$
where on the last step we use the following consequences of~\thetag{4}:
$$
\sum_{l=0}^nA_l(0)H_n
=\sum_{l=0}^nA_l(0)H_{2n}
=0.
$$
Since
$$
A_l(0)\cdot\biggl(\frac1{\frac n2-l}-6H_{n-l}+6H_l
-2H_{n+l}+2H_{2n-l}\biggr)
=-\frac{\d}{\d l}A_l(0),
$$
after developing all Pochhammer symbols in the $\epsilon\to0$
form of~\thetag{3} we arrive at the identity from
\cite{KR, Section~13}:
$$
\align
&
-\sum_{l=0}^n\frac{\d}{\d l}
\Bigl(\frac n2-l\Bigr)\binom nl^4\binom{n+l}n^2\binom{2n-l}n^2
\\ &\quad
=(-1)^n\sum_i\binom ni^2\binom{2n-i}n
\sum_j\binom{n+j-i}n\binom nj^2\binom{n+j}n.
\tag5
\endalign
$$

Clearly, the left-hand side of Andrews's formula is symmetric
with respect to the group of parameters $b_1,c_1,b_2,c_2,b_3,c_3$.
Therefore, setting as before $a=-n-2\epsilon$, $m=n$, and all
the parameters of the group to be $-n-\epsilon$, except the following two:
\roster
\item"(a)" $b_1=c_1=n-\epsilon+1$;
\item"(b)" $b_2=c_2=n-\epsilon+1$;
\item"(c)" $b_3=c_3=n-\epsilon+1$;
\item"(d)" $c_1=c_2=n-\epsilon+1$;
\item"(e)" $c_2=c_3=n-\epsilon+1$;
\item"(f)" $c_1=c_3=n-\epsilon+1$
\endroster
(the last case corresponds to the above identity~\thetag{5}),
we arrive at the five more representations of the left-hand side
of~\thetag{5}:
$$
\allowdisplaybreaks
\align
(-1)^nu_n
&=-\sum_{l=0}^n\frac{\d}{\d l}
\Bigl(\frac n2-l\Bigr)\binom nl^4\binom{n+l}n^2\binom{2n-l}n^2
\\
&=(-1)^n\sum_i(-1)^i\binom{3n+1}i\binom{2n-i}n^2
\sum_j\binom{n+j-i}n\binom nj^2\binom{2n-j}n
\\
&=(-1)^n\sum_i(-1)^i\binom{n+i}n^3
\sum_j(-1)^j\binom{3n+1}{j-i}\binom{2n-j}n^3
\\
&=\sum_i\binom ni^2\binom{n+i}n
\sum_j(-1)^j\binom{n+j-i}n\binom{n+j}n^2\binom{3n+1}{n-j}
\\
&=(-1)^n\sum_i\binom ni\binom{n+i}n\binom{2n-i}n
\sum_j\binom n{j-i}\binom nj\binom{2n-j}n^2
\\
&=(-1)^n\sum_i\binom ni\binom{n+i}n^2
\sum_j\binom n{j-i}\binom nj\binom{n+j}n\binom{2n-j}n.
\endalign
$$


\enddocument